\documentclass{article}
\usepackage{amsmath,amsthm,amssymb}
\usepackage[authoryear,math]{artratex}
\title{On the log-hyperconvexity index of pseudoconvex H\"{o}lder domains in $\mathbb{C}^n$ }
\author{Tianlong Yu }

\begin{document}
\maketitle%

\textbf{Abstract}:
In this note we prove that every bounded pseudoconvex domain in $\mathbb{C}^n$ with H\"{o}lder boundary admits a positive log-hyperconvexity index. We also get a better estimate for the plurisubharmonic exhaustion function.

\section{Introduction}
A bounded pseudoconvex domain $\Omega$ in $\mathbb{C}^n$ is called hyperconvex if there exists $0>\phi \in PSH(\Omega)$ such that for every $c<0$, $\phi^{-1}([-\infty, c))$ is relative compact in $\Omega$. 
In their pioneer work [7], Diederich and Fornaess showed that every bounded pseudoconvex domain in $\mathbb{C}^n$ with $C^2$ boundary is hyperconvex. Actually they proved that there exists a smooth psh function which is comparable to -$d^{\alpha}$, for some $1>\alpha>0$, where $d$ is the boundary distance. In [9] Kerzman and Rosay extended the above result to show that every bounded pseudoconvex domain in $\mathbb{C}^n$ with $C^1$ boundary is hyperconvex. In [6] Demailly proved that every bounded pseudoconvex domain in $\mathbb{C}^n$ with Lipschitz boundary is hyperconvex. Actually he constructed a smooth psh function which is comparable to $(\log d)^{-1}$. In [8] Harrington further generalized Demailly's result to show that every bounded pseudoconvex domain in $\mathbb{C}^n$ with Lipschitz boundary also admits a psh function which is comparable to -$d^{\alpha}$, for some $1>\alpha>0$. In [2], by combining the localization principle for hyperconvex domains in $\mathbb{C}^n$ proved in [9] and a new estimate of relative extremal function, Chen showed that every bounded pseudoconvex domain in $\mathbb{C}^n$ with H\"{o}lder boundary is hyperconvex. Actually he proved that $\textit{locally}$ there exists a psh function(the local relative extremal function) which has a lower bound comparable to $-(-\log d)^{-\tau}$ , for some $1>\tau>0$. In [3], by combining the results in [2] and a new psh Hopf lemma and the patching methods proved in [5] , Chen and Xiong obtained a lower bound of the global relative extremal function in the H\"{o}lder boundary case.

In [4] , as the hyperconvexity index introduced by Chen[1], Chen and Zheng introduced the concept of log-hyperconvexity index as follows: 
\begin{definition}
For a bounded domain $\Omega$ in $\mathbb{C}^n$, define
\begin{equation}
 \alpha_{l}(\Omega)=\sup\{\alpha \geq 0; \exists  0>\phi \in C(\Omega)\cap PSH(\Omega) s.t.-\phi \leq C(-\log d)^{-\alpha}\}.   
\end{equation}
 as the log-hyperconvexity index of $\Omega$.
 \end{definition}
 
 Therefore the main results in [2] read as every bounded pseudoconvex domain in $\mathbb{C}^n$ with H\"{o}lder boundary $\textit{locally}$ has a positive log-hyperconvexity index. And in [2][3][4], Chen raised the question that whether every bounded pseudoconvex domain in $\mathbb{C}^n$ with H\"{o}lder boundary has a ($\textit{global}$) positive log-hyperconvexity index. In this note, we give a positive answer to this question by constructing a plurisubharmonic exhaustion function with better growth estimate near the boundary:
 
\begin{theorem}
 Let $\Omega$ be a bounded pseudoconvex domain in $\mathbb{C}^n$ with H\"{o}lder boundary. Then for some $\tau>0$, there exists a continuous psh function $w(z)$ on $\Omega$ and $M_1, M_2>0$ such that 
 \begin{equation}
 -\frac{M_1}{(-\log d(z))^{\tau}} \frac{1}{\log(-\log d(z))}\leq w(z)\leq -\frac{M_2}{(-\log d(z))^{\tau}} \frac{1}{\log(-\log d(z))}.    
 \end{equation}
 
\end{theorem}
 
 The proof of the theorem is a standard psh patching argument as used by Kerzman-Rosay[9], Demailly[6], Coltiou-Mihalche[5], Harrington[8]. Here we work twice the patching process, once locally and once globally.
\section{Proof of Theorem 1}

 We say $\Omega$ has H\"{o}lder boundary or $\Omega$ is a H\"{o}lder domain if locally $\partial \Omega$ can be written as a H\"{o}lder continuous graph under certain holomorphic coordinates.
Hence we can choose a finite collection of points $p_j\in \partial \Omega, r_j>0, U_j''\subset\subset U_j'\subset\subset U_j$ open neighborhoods of $p_j$ such that $\partial \Omega \subseteq \cup U_j''$ and the following holds after local holomorphic coordinate changes:
\begin{equation}\label{eq:2.1}
U_j=B(0,10r_j);
\end{equation}
\begin{equation}\label{eq:2.2}
 U_j'=B(0,2r_j);   
\end{equation}
\begin{equation}\label{eq:2.3}
 U_j\cap \Omega=\{z\in B(0,10r_j)|\mathrm{Im}z_n<g_j(z',\mathrm{Re}z_n)\}  . 
\end{equation}
\noindent where $z'=(z_1,\cdots,z_{n-1})$ and $g_j$ is a H\"{o}lder continuous function defined on $\{(z', x_n)\in \mathbb{C}^{n-1}\times \mathbb{R}| |z'|^2+x_{n}^2<100r_{j}^2\}$ with a H\"{o}lder index $\beta_j \in (0,1)$ and satisfies $g_j(0', 0)=0$.(Here we assume $p_j=0$ in local coordinates) For $t>0$, let 
\begin{equation}\label{eq:2.4}
\Omega_{t,j}=\{z\in B(0,10r_j)|\mathrm{Im}z_n<g_j(z',\mathrm{Re}z_n)+t\}.   
\end{equation}

\noindent which are also pseudoconvex domains containing $U_j\cap \Omega =:\Omega_{0,j}$. Let $d_{t,j}$ be the boundary distance of $\Omega_{t,j}, C_{t,j}=d(\partial \Omega_{0,j}\cap B(0,8r_j),\partial \Omega_{t,j}\cap B(0,9r_j))$. Then by (\ref{eq:2.3})and(\ref{eq:2.4}), as in the proof of Theorem1.1 in [2], we can find $\gamma_j>1$ which depends on $\beta_j$, such that 
\begin{equation}\label{eq:2.5}
    t^{\gamma_j}\leq C_{t,j}\leq t.
\end{equation}
\noindent for all small $t$.We may assume $\gamma_j=\gamma$ for all $j$. In the following we first work locally in $B(0,2r_j)$, so we omit $j$ for the sake of convenience. So $\Omega_t=\Omega_{t,j}; r=r_j; d_t=d_{t,j}; C_t=C_{t,j}$. For $z\in \Omega_0$,let 
\begin{equation}\label{eq:2.6}
v_n(z)=\frac{-\log d_{t_n}(z)-(1+\epsilon_{0})(-\log C_{t_n})}{a_{n}(-\log C_{t_n})};
\end{equation}
\begin{equation}\label{eq:2.7}
E_n=\{z\in \Omega_0|t_{n+1}'\leq d_0(z)<t_{n}'\};    
\end{equation}

\noindent where $a_{n}, \epsilon_0, t_{n}, t_{n}'$ are positive constants which will be determined later. For any positive integer n and $z\in E_n$, we define 
\begin{equation}\label{eq:2.8}
  v(z)=\max\{v_{1}(z),\cdots, v_{n}(z)\}.  
\end{equation}

\noindent We want to choose $t_n, t_{n}', a_n, \epsilon_0$ properly such that $v(z)$ is a well-defined psh function on $\Omega_0$. For that purpose, we only need to check $v_{n}(z)>v_{n+1}(z)$, when $d_0(z)=t_{n+1}'$.

We claim that for all $z\in \Omega_0\cap B(0,2r)$, and $t>0$ sufficiently small
\begin{equation}\label{eq:2.9}
  d_t(z)\geq d_0(z)+C_t;  
\end{equation}
\begin{equation}\label{eq:2.10}
  d_t(z)\leq d_0(z)+t.  
\end{equation}
\noindent Indeed, let $z_j\in \partial \Omega_{t}\cap B(0,10r)$ such that lim$_{j\to \infty}|z-z_j|=d_t(z)$, then there exists $z_j'\in \partial \Omega_{0}\cap B(0,10r)$ such that $|z-z_j|=|z-z_j'|+|z_j'-z_j|$, when $t$ small enough, we can achieve that all $z_j$ and $z_j'$ are in $B(0, 8r)$, thus $|z_j'-z_j|\geq C_t$, hence (\ref{eq:2.9}) holds. For (\ref{eq:2.10}), let $z'\in \partial \Omega_{0}\cap B(0,10r)$ such that $d_0(z)=|z-z'|$, let $z'_t=z'+(0',it)$, then for small enough $t$, $z'_t\in \partial \Omega_{t}\cap B(0,10r)$, thus $d_t(z)\leq |z-z'_t|\leq |z-z'|+t=d_0(z)+t$.

By (\ref{eq:2.5}),(\ref{eq:2.6}),(\ref{eq:2.9})and(\ref{eq:2.10}), we have 
\begin{equation}\label{eq:2.11}
 v_{n+1}(z)\leq -\frac{1}{a_{n+1}}+\frac{\log(t_{n+1}'+C_{t_{n+1}})}{a_{n+1}\log C_{t_{n+1}}}\leq -\frac{1}{a_{n+1}}+\frac{\log(t_{n+1}'+t_{n+1})}{a_{n+1}\log t_{n+1}};
\end{equation}
\begin{equation}\label{eq:2.12}
 v_{n}(z)\geq -\frac{1+\epsilon_0}{a_{n}}+\frac{\log(t_{n+1}'+t_{n})}{a_{n}\log C_{t_{n}}}\geq -\frac{1+\epsilon_0}{a_{n}}+\frac{\log(t_{n+1}'+t_{n})}{a_{n}\log(t_{n}^{\gamma})}.   
\end{equation}

\noindent when $z\in \Omega_0\cap B(0,2r)$, and $d_0(z)=t_{n+1}'$. Here we use the monotonicity of the function $\frac{\log(x+c)}{\log x}$ for $c>0$, $x>0$ and $x+c<e^{-1}$. Now we fix $0<t_1<1$. Let 
\begin{equation}\label{eq:2.13}
 t_{n+1}=t_{n}^{2\gamma}=\cdots=t_{1}^{(2\gamma)^n};   
\end{equation}
\begin{equation}\label{eq:2.14}
t_{n+1}'=t_n.    
\end{equation}

\noindent Then by (\ref{eq:2.13})and(\ref{eq:2.14}), we have 
\begin{equation}\label{eq:2.15}
  \frac{\log(t_{n+1}'+t_{n+1})}{\log t_{n+1}}\leq \frac{\log((t_1)^{(2\gamma)^{n-1}})}{\log t_{n+1}}=\frac{1}{2\gamma}<1;  
\end{equation}
\begin{equation}\label{eq:2.16}
  \frac{\log(t_{n+1}'+t_{n})}{\log(t_{n}^{\gamma})}\geq \frac{\log2}{(2\gamma)^{n-1}\gamma(\log t_1)}+\frac{1}{\gamma}.  
\end{equation}

\noindent We choose $n$ large and $\epsilon_0$ small such that 
\begin{equation}\label{eq:2.17}
  1-\frac{1}{2\gamma}>1+\epsilon_0-\frac{1}{\gamma}-\frac{\log2}{(2\gamma)^{n-1}\gamma(\log t_1)}.  
\end{equation}

\noindent By (\ref{eq:2.11}),(\ref{eq:2.12}),(\ref{eq:2.15}),(\ref{eq:2.16})and(\ref{eq:2.17}), we can find $\lambda>1$, such that for $a_{n}=\lambda^{n}$, we have $v_{n}(z)>v_{n+1}(z)$, when $z\in \Omega_0\cap B(0,2r)$, and $d_0(z)=t_{n+1}'$.

By (\ref{eq:2.9}), for $z\in \Omega_0\cap B(0,2r)$, $d_t(z)\geq C_t$ for any $t>0$. Hence, by (\ref{eq:2.6}),(\ref{eq:2.7})and(\ref{eq:2.8}), for $z\in E_n\cap B(0,2r)$, we have
\begin{equation}\label{eq:2.18}
 v(z)\leq -\frac{\epsilon_0}{\lambda^{n}}.   
\end{equation}

\noindent On the other hand, for $z\in E_n\cap B(0,2r)$
\begin{equation}\label{eq:2.19}
 v(z)\geq v_{n}(z)\geq -\frac{1+\epsilon_0}{\lambda^{n}};   
\end{equation}
\begin{equation}\label{eq:2.20}
 (2\gamma)^{n-1}\mathrm{log} t_1\leq\mathrm{log}d_0(z)<(2\gamma)^{n-2}\mathrm{log}t_1.  
\end{equation}
\noindent Thus by (\ref{eq:2.18}),(\ref{eq:2.19})and(\ref{eq:2.20}), there exists $C_1$,$C_2>0$, $\tau'>0$ such that for all $z\in \Omega_0\cap B(0,2r)$
\begin{equation}\label{eq:2.21}
  -\frac{C_2}{(-\log d_0(z))^{\tau'}}
 \leq v(z)\leq -\frac{C_1}{(-\log d_0(z))^{\tau'}}.  
\end{equation}
 
\noindent Now we again work locally for fixed $j$, but with positive $t$. In other words, we replace the above $\Omega_0$ by $\Omega_{t,j}$ which are also H\"{o}lder domains with the same H\"{o}lder index as $\Omega_{0,j}$. Therefore, after some rescaling, we get, by the above arguments, a family of psh functions $v_{t,j}$ defined on $\Omega_{0,j}\cap B(0,2r_j)$ for all small positive $t$, such that there exists $C>1$, $\tau_{j}>0$, independant of $t$, with the following estimates hold on $\Omega_{0,j}\cap B(0,2r_j)$
\begin{equation}\label{eq:2.22}
  -\frac{C}{(-\log d_{t,j}(z))^{\tau_j}}\leq v_{t,j}(z)\leq -\frac{1}{(-\log d_{t,j}(z))^{\tau_j}}.  
\end{equation}

\noindent Let $\tau_{j_0}=\min\{\tau_j\}$, and $w_{t,j}(z)=-(-v_{t,j}(z))^{\frac{\tau_{j_0}}{\tau_j}}$, then $w_{t,j}(z)$ are psh functions and
\begin{equation}\label{eq:2.23}
  -\frac{C^{\frac{\tau_{j_0}}{\tau_j}}}{(-\log d_{t,j}(z))^{\tau_{j_0}}}\leq w_{t,j}(z)\leq -\frac{1}{(-\log d_{t,j}(z))^{\tau_{j_0}}}.  
\end{equation}

\noindent Hence 
\begin{equation}\label{eq:2.24}
    -\frac{\tau_{j_0}}{\tau_j}\log C+\tau_{j_0}\log(-\log d_{t,j}(z))\leq -\log(-w_{t,j}(z))\leq \tau_{j_0}\log(-\log d_{t,j}(z)).
\end{equation}

\noindent By (\ref{eq:2.5}),(\ref{eq:2.9})and(\ref{eq:2.10}), for $z\in \Omega_{0,j}\cap B(0,2r_j)$, we have
\begin{equation}\label{eq:2.25}
  d_{0,j}(z)+t^{\gamma}\leq d_{t,j}(z)\leq d_{0,j}(z)+t . 
\end{equation}
 
 \noindent Thus there exists $C_3>0$, independant of $t, j, z$, such that 
 \begin{equation}\label{eq:2.26}
  0\leq \log(-\log d_{t,j}(z))-\log(-\log(d(z)+t))\leq C_3.  
 \end{equation}
 
 \noindent (Here we actually measure the boundary distance in different local holomorphic coordinates, but they are comparable, so we can use the standard Euclidean distance $d(z)$.)

Now, we use Richberg's patching technique[10]. Let
\begin{equation}\label{eq:2.27}
   w_t(z)=\max_{\{j:z\in U_j'\}}\{-\log(-w_{t,j}(z))+\eta_j(z)+M|z|^2\}. 
\end{equation}

\noindent where $\eta_j(z)>\tau_{j_0}C_3+\frac{\tau_{j_0}}{\tau_j}\log C$ when $z\in U_j''$, with compact support in $U_j'$, $M$ is large enough to make $w_t$ psh.
By (\ref{eq:2.2}),(\ref{eq:2.24})and(\ref{eq:2.26}), $w_t$ is a well-defined psh function on $\Omega$ (near the boundary) and, after adding some constants, there exists $C_4>0$ such that
\begin{equation}\label{eq:2.28}
  -C_4+\tau_{j_0}\log(-\log(d(z)+t))\leq w_t(z)\leq \tau_{j_0}\log(-\log(d(z)+t^{\gamma})).  
\end{equation}

\noindent Next we do the psh patching process globally on $\Omega$. For $z\in\Omega$, let
\begin{equation}\label{eq:2.29}
w_n(z)=\frac{w_{t_n}(z)-\tau_{j_0}\log(-\log(t_{n}^{\gamma}))-\epsilon_1}{\tau_{j_0}\log(-\log(t_{n}^{\gamma}))b_n};
\end{equation}
\begin{equation}\label{eq:2.30}
E_{n}'=\{z\in \Omega|t_{n+1}''\leq d(z)<t_{n}''\};
\end{equation}
\noindent where $\epsilon_1>0$ is a fixed constant, and $t_n, t_{n}'', b_n$ are positive constants which will be determined later. For any positive integer n and $z\in E_{n}'$, we define
\begin{equation}\label{eq:2.31}
  w(z)=\max\{w_{1}(z),\cdots, w_{n}(z)\}.  
\end{equation}
\noindent Again we shall choose $t_{n}, t_{n}'', b_n$ properly to make $w(z)$ a well-defined psh function on $\Omega$. For that purpose, we only need to check $w_{n}(z)>w_{n+1}(z)$, when $d(z)=t_{n+1}''$. By (\ref{eq:2.27}),(\ref{eq:2.28})and(\ref{eq:2.29}), we have 
\begin{equation}\label{eq:2.32}
   w_n(z)\geq -\frac{1}{b_n}[1-\frac{\tau_{j_0}\log(-\log(d(z)+t_n))-\epsilon_1-C_4}{\tau_{j_0}\log(-\log(t_{n}^{\gamma}))}]; 
\end{equation}
\begin{equation}\label{eq:2.33}
  w_{n+1}(z)\leq -\frac{1}{b_{n+1}}[1-\frac{\log(-\log(d(z)+t_{n+1}^{\gamma}))}{\log(-\log(t_{n+1}^{\gamma}))}].  
\end{equation}

\noindent Let $k$ be a postive integer, $0<t_1<1$ be fixed and let
\begin{equation}\label{eq:2.34}
 t_{n+1}''=t_{n+1}^{\frac{1}{k}}-t_{n+1}^{\gamma};  
\end{equation}
\begin{equation}\label{eq:2.35}
 t_{n+1}=t_{n}^k=\cdots=t_{1}^{k^n}.   
\end{equation}

\noindent Then $\lim_{n\to \infty}t_{n}''=0$, and for large $k$ and $n$, we have $t_{n+1}''<t_{n}''$. When $d(z)=t_{n+1}''$, since $t_{n+1}''<t_{n}$, we have, by (\ref{eq:2.32})and(\ref{eq:2.33}),
\begin{equation}\label{eq:2.36}
1-\frac{\tau_{j_0}\log(-\log(d(z)+t_n))-\epsilon_1-C_4}{\tau_{j_0}\log(-\log(t_{n}^{\gamma}))}<\frac{\log(\frac{\log t_{n}^{\gamma}}{\log(2t_n)})}{\log(-\log t_{n}^{\gamma})}+\frac{\epsilon_1+C_4}{\tau_{j_0}\log(-\log(t_{n}^{\gamma}))};    
\end{equation}
\begin{equation}\label{eq:2.37}
1-\frac{\log(-\log(d(z)+t_{n+1}^{\gamma}))}{\log(-\log(t_{n+1}^{\gamma}))}=\frac{\log(k\gamma)}{\log(-\log t_{n+1}^{\gamma})}=\frac{\log(k\gamma)}{\log(-\log t_{n}^{\gamma})+\log k}.    
\end{equation}

\noindent By (\ref{eq:2.32}),(\ref{eq:2.33}),(\ref{eq:2.36})and(\ref{eq:2.37}), we see that if we choose $k$ such that log$k>\frac{\epsilon_1+C_4}{\tau_{j_0}}$, there exists $\lambda'>1$ such that, by letting $b_{n}=(\lambda')^n$, we have $w_{n}(z)>w_{n+1}(z)$ for large $n$, when $d(z)=t_{n+1}''$. By (\ref{eq:2.28}),(\ref{eq:2.29}),(\ref{eq:2.30})and(\ref{eq:2.31}), there exists a   positive constant $C_5$ independant of $n$ such that for $z\in E_{n}'$, we have
\begin{equation}\label{eq:2.38}
 w(z)\leq -\frac{\epsilon_1}{(\lambda')^{n}\tau_{j_0}\log(-\log(t_{n}^{\gamma}))};   
\end{equation}
\begin{equation}\label{eq:2.39}
 w(z)\geq w_n(z)\geq -\frac{C_5}{(\lambda')^{n}\tau_{j_0}\log(-\log(t_{n}^{\gamma}))};   
\end{equation}
\begin{equation}\label{eq:2.40}
  \log t_{n+1}''\leq \log d(z)< \log t_{n}''.
\end{equation}

\noindent Therefore, by (\ref{eq:2.34}),(\ref{eq:2.35}),(\ref{eq:2.38}),(\ref{eq:2.39})and(\ref{eq:2.40}), there exists $M_1, M_2>0, \tau>0$ such that for all $z\in \Omega$
\begin{equation}\label{eq:2.41}
 -\frac{M_1}{(-\log d(z))^{\tau}} \frac{1}{\log(-\log d(z))}\leq w(z)\leq -\frac{M_2}{(-\log d(z))^{\tau}} \frac{1}{\log(-\log d(z))}.   
\end{equation}

\noindent As all the boundary distance functions are continuous and the whole patching process only involves maximum of a finite collection of functions, $w(z)$ is continuous psh and therefore we complete the proof.

\section*{Acknowledgement}
I would like to thank my supervisor professor Xiangyu Zhou who encouraged me to study problems related to hyperconvexity. I also thank Zhiqiang Wu and Hui Yang drawing my attention to related papers.

\section*{Declarations}
\textbf{Ethical Approval}: 
\\
This declaration is not applicable.
\\
\textbf{Funding}:
\\
The author has no funding support.
\\
\textbf{Availability of data and materials}:
\\
This declaration is not applicable.

Tianlong Yu: School of Mathematical Sciences, Peking University, Beijing, 100871, China.

National Center for Mathematics and Interdisciplinary Sciences, CAS, Beijing, 100190, China.

E-mail address: yutianlong940309@163.com

\end{document}